\documentclass[11pt,a4paper]{article}
\usepackage[margin=2.5cm]{geometry}
\usepackage[T1]{fontenc}
\usepackage{lmodern}
\usepackage{amsmath,amssymb,amsthm,mathtools}
\usepackage{thmtools}
\usepackage{microtype}
\usepackage{enumitem}
\usepackage{xcolor}

\usepackage{float,tabularx,array}
\floatstyle{ruled}
\newfloat{algorithm}{tbp}{loa}
\floatname{algorithm}{Algorithm}
\newcolumntype{L}{>{\raggedright\arraybackslash}X}
\newcommand{\algline}[2]{#1 & #2\\[1pt]}
\usepackage{hyperref}
\usepackage[nameinlink,capitalize]{cleveref}
\crefname{algorithm}{Algorithm}{Algorithms}

\hypersetup{
  colorlinks=true,
  linkcolor=blue!50!black,
  citecolor=blue!50!black,
  urlcolor=blue!50!black,
  pdftitle={Every r-regular Hamiltonian graph with r>16 contains a second Hamilton cycle},
  pdfauthor={Haiyang Liu, Xu Liu, Bo Ning}
}

\allowdisplaybreaks[1]
\newtheorem{theorem}{Theorem}[section]
\newtheorem{lemma}[theorem]{Lemma}

\newtheorem{corollary}[theorem]{Corollary}
\newtheorem{claim}{Claim}
\theoremstyle{definition}

\newtheorem{conjecture}[theorem]{Conjecture}
\numberwithin{equation}{section}

\makeatletter
\renewenvironment{proof}[1][\proofname]{%
  \par\pushQED{\qed}\normalfont
  \topsep6\p@\@plus6\p@\relax
  \trivlist
  \item[\hskip\labelsep\bfseries #1\@addpunct{.}]\ignorespaces
}{%
  \popQED\endtrivlist\@endpefalse
}
\makeatother

\newcommand{\PrOmega}{\Pr\nolimits_{\Omega}}
\newcommand{\PrMT}{\Pr\nolimits_{\mathrm{MT}}}

\title{\fontsize{14}{18}\selectfont
Every $r$-regular Hamiltonian graph with $r\geq 17$\\
contains a second Hamilton cycle}
\author{
    Haiyang Liu\thanks{College of Cryptology and Cyber Security, Nankai University, Tianjin 300350, P.R. China.}
    \and
    Xu Liu\thanks{College of Computer Science, Nankai University, Tianjin 300350, P.R. China.}
    \and 
    Bo Ning\thanks{Corresponding author. College of Computer Science \& College of Cryptology and Cyber Security, Nankai University, Tianjin 300350, P.R. China. E-mail: \texttt{bo.ning@nankai.edu.cn} (B. Ning). Partially supported by the National Natural Science Foundation of China (No. 12371350). \href{https://orcid.org/0000-0002-9622-5567}{ORCID: 0000-0002-9622-5567}.}}
\date{}

\begin{document}
\maketitle

\begin{abstract}
A conjecture of Sheehan states that every Hamiltonian $4$-regular graph contains a second Hamilton cycle. A natural problem is to find the smallest positive integer $k$ such that every Hamiltonian $r$-regular graph contains a second Hamilton cycle when $r\geq k$. Thomassen proved that $k\leq 300$. Improving Thomassen's result, Haxell, Seamone, and Verstra{\"e}te proved that $k\leq 23$. We improve this bound by proving that every $r$-regular Hamiltonian graph contains a second Hamilton cycle whenever $r\geq 17$. Our new idea is to establish a transversal theorem for set systems on a cycle using Harris's results on the Moser--Tardos algorithm for Local Lemma. We then apply Thomassen's theorem on independent dominating sets to obtain a second Hamilton cycle.
For $r\ge 18$, we also give some algorithms to find such a transversal.
\end{abstract}

\section{Introduction}

Throughout the paper, all graphs are finite, simple, and undirected. A graph is \emph{uniquely Hamiltonian} if it contains exactly one Hamilton cycle. A classical question asks whether an $r$-regular graph can be uniquely Hamiltonian when $r\ge3$. The restriction on $r$ is necessary, since every cycle is uniquely Hamiltonian.

A theorem of C.~A.~B. Smith, first published by Tutte~\cite{Tutte1946} in 1946, states that every edge of a cubic graph lies in an even number of Hamilton cycles. In particular, every Hamiltonian cubic graph contains a second Hamilton cycle. Thomason~\cite{Thomason1978} extended this parity theorem to graphs in which every vertex has odd degree. That is, every Hamiltonian $r$-regular graph with an odd $r\ge3$ contains a second Hamilton cycle. For even degrees, Sheehan~\cite{Sheehan1975} proposed the following conjecture in 1975.

\begin{conjecture}[\cite{Sheehan1975}]
Every Hamiltonian $4$-regular graph contains a second Hamilton cycle.
\end{conjecture}

Thomassen~\cite{Thomassen1998} proved that every Hamiltonian $r$-regular graph with $r\ge300$ contains a second Hamilton cycle. His proof combines a sufficient condition from his earlier work~\cite{Thomassen1997} with the Lov\'asz Local Lemma. Using the Lov\'asz Local Lemma, Hor\'ak and Stacho \cite{HorakStacho2000} improved Thomassen's bound by showing that every $r$-regular Hamiltonian graph with $r\ge 73$ contains a second Hamiltonian cycle. 
Haxell, Seamone, and Verstra{\"e}te~\cite{HaxellSeamoneVerstraete2007} subsequently improved Thommassen's bound for regular graphs by using a different random construction. Their main result is the following.

\begin{theorem}[{\cite{HaxellSeamoneVerstraete2007}, Theorem~1}]
Every Hamiltonian $r$-regular graph with $r\geq 23$ contains a second Hamilton cycle.
\end{theorem}

We improve this degree bound as follows.

\begin{theorem}\label{thm:main}
Let $r\geq 17$. Every $r$-regular graph containing a Hamilton cycle contains a second Hamilton cycle.
\end{theorem}

The case $r=17$ follows from Thomason's theorem. Our argument proves the result for every $r\ge18$ by establishing a more general theorem on transversals of a family of sets on a cycle.

Fix a Hamilton cycle $C$ in a graph $G$, and color the edges of $C$ red and the remaining edges green. A set $S\subseteq V(G)$ is \emph{red-stable} if no two vertices of $S$ are joined by a red edge, and \emph{green-dominating} if every vertex of $V(G)\setminus S$ has a green neighbor in $S$. For a family $\mathcal F$ of subsets of $V(C)$, a \emph{red-stable transversal} of $\mathcal F$ is a set that is independent in $C$ and meets every member of $\mathcal F$.

\begin{theorem}\label{thm:transversal}
Let $G$ be a graph with a fixed Hamilton cycle $C$.
Let $k\ge17$, and let $\mathcal D$ be a family of $k$-element subsets of $V(C)$.
Color the edges of $C$ red and the remaining edges of $G$ green.
If every vertex of $C$ belongs to at most $k$ members of $\mathcal D$, then $\mathcal D$ has a red-stable transversal.
\end{theorem}

Red-stable green-dominating sets need not exist in small regular degrees. Goedgebeur, Jooken, Lo, Seamone, and Zamfirescu~\cite{GoedgebeurEtAl2024} showed that, for each $r\in\{5,6\}$, there are infinitely many $r$-regular graphs with a prescribed Hamilton cycle $C$ for which no such set exists.

For a graph $G$, we write $V(G)$ and $E(G)$ for its vertex set and edge set, respectively, and $d_G(v)$, $N_G(v)$, and $N_G[v]$ for the degree, neighborhood, and closed neighborhood of a vertex $v$. For $F\subseteq E(G)$, the graph $G-F$ is obtained by deleting the edges in $F$ and retaining all vertices. In particular, $H=G-E(C)$ denotes the spanning subgraph consisting of the green edges.

The paper is organized as follows. In \Cref{sec:lll}, we recall the Local Lemma and the resampling results needed in the proof. In \Cref{sec:transversal}, we give the structural reduction, establish the weighted path estimates, prove \cref{thm:transversal}, and deduce \cref{thm:main}.
In \Cref{sec:algorithm}, we present algorithms for constructing such a transversal and analyze their computational complexity.

\subsection{Outline of proofs of main results}
\label{sec:outline}

We first explain how \cref{thm:transversal} implies \cref{thm:main}, and then outline the proof of the transversal theorem. The structural step is a criterion of Thomassen~\cite{Thomassen1997,Thomassen1998}: if a graph $G$ with a fixed Hamilton cycle $C$ has a set that is independent in $C$ and dominating in $G-E(C)$, then $G$ contains a Hamilton cycle distinct from $C$. Such a set is often called a \emph{$C$-independent dominating set}; in our terminology, it is red-stable and green-dominating.

For our application, fix a Hamilton cycle $C$ in an $r$-regular graph $G$ with $r\ge18$, and let $H=G-E(C)$ and $k=r-1$. Since $H$ is $(r-2)$-regular, every closed neighborhood $N_H[v]$ has size $k$. Moreover, a vertex $y$ belongs to $N_H[v]$ precisely when $v\in N_H[y]$, so every vertex belongs to at most $k$ distinct closed neighborhoods. By \cref{thm:transversal}, these neighborhoods have a red-stable transversal $S$. Every vertex outside $S$ then has a neighbor in $S$ through an edge of $H$. Thus $S$ is independent in $C$ and dominating in $H$, and Thomassen's criterion gives a second Hamilton cycle. This reduction is recorded in \cref{lem:reduction}.

To prove \cref{thm:transversal}, we construct a random set that is independent in $C$ for every outcome. Delete a vertex $\rho$ from $C$, and choose a bipartition $V(P)=A\mathbin{\dot\cup}B$ of the resulting path $P=C-\rho$ so that at most one endvertex of $P$ lies in $A$. Select each vertex of $B$ independently with probability $p=13/40$, let $T$ be the set of selected vertices, and define
$S=T\cup\{a\in A:N_P(a)\cap T=\varnothing\}$.
The set $S$ is independent in $P$, and hence in $C$ because $\rho\notin S$, as verified in \cref{lem:independent}. It remains to arrange that $S$ meets every member of $\mathcal D$.

For $D\in\mathcal D$, let
$Z_D=D\cap B,
\mathcal F_D=\{N_P(a):a\in D\cap A\}$.
The event $S\cap D=\varnothing$ occurs precisely when no vertex of $Z_D$ is selected and $T$ meets every member of $\mathcal F_D$. In \cref{lem:encoding}, we express this event as a union of atomic events indexed by the minimal transversals $R$ of $\mathcal F_D$ that avoid $Z_D$. The corresponding atomic event prescribes value $0$ on $Z_D$ and value $1$ on $R$.

The main estimate uses the path structure of $\mathcal F_D$. Order the vertices of $B$ along $P$ and join consecutive vertices to form an auxiliary path $Q$. Every member of $\mathcal F_D$ is the endvertex set of an edge of $Q$, with the possible exception of one singleton at an endvertex. For $\alpha=29/40$ and $\beta=14/40$, \cref{lem:path-estimate} bounds the sum of $\beta^{\lvert R\rvert}$ over all minimal transversals of such a family by $\alpha^{\lvert\mathcal F_D\rvert}$. It also bounds the corresponding sum over transversals containing a prescribed vertex. The first bound controls the total weight of the atomic events arising from one set $D$; the second controls the weight of events requiring a fixed variable to be $1$.

We separate the sets avoiding $\rho$ from those containing it, writing
$\mathcal D_0=\{D\in\mathcal D:\rho\notin D\},
\mathcal D_1=\{D\in\mathcal D:\rho\in D\}$.
For the atomic events arising from $\mathcal D_0$, we use the weights $\mu(E(D,R))=\alpha^{\lvert Z_D\rvert}\beta^{\lvert R\rvert}$. The weighted path estimates, together with the assumption that each vertex belongs to at most $k$ members of $\mathcal D$, give bounds on the total weights of events imposing either value on a fixed variable. These bounds are proved in \cref{lem:loads} and used in \cref{lem:resampling-check} to verify Harris's criterion~\cite{Harris2016} for the Moser--Tardos algorithm~\cite{MoserTardos2010}. The required inequalities hold at $k=17$, and the relevant bounds decrease as $k$ increases. Consequently, the algorithm terminates with probability one, and every terminal assignment gives a set $S$ meeting all members of $\mathcal D_0$.

Finally, the occurrence bound gives $\lvert\mathcal D_1\rvert\le k$. We use Harris's estimate for the terminal distribution~\cite{Harris2016} to show that each $D\in\mathcal D_1$ is missed with probability at most $\alpha^{k-1}$. A union bound gives
$\PrMT\bigl(\text{$S$ misses some member of $\mathcal D_1$}\bigr)
\le k\alpha^{k-1}<1$,
where $\PrMT$ denotes probability under the terminal distribution. Thus, some terminal assignment meets every member of $\mathcal D_1$ as well, completing the proof of \cref{thm:transversal}. The use of the terminal distribution allows us to control these remaining events after the resampling algorithm has enforced the other constraints; see also Harris and Srinivasan~\cite{HarrisSrinivasan2017} for a systematic study of this distribution.

\section{Local Lemma and resampling tools}\label{sec:lll}

Throughout this section, $[n]=\{1,\ldots,n\}$. For an event $A$, $\overline{A}$ denotes the complement of $A$. In 1975, Erd\H{o}s and Lov\'asz~\cite{ErdosLovasz1975} introduced the Lov\'asz Local Lemma, and Spencer~\cite{Spencer1977} gave the following symmetric and asymmetric forms in 1977.

Let $I$ be a finite index set, and let $A_i$ be an event for each $i\in I$. A graph $D$ with vertex set $I$ is a \emph{dependency graph} for $\{A_i:i\in I\}$ if each $A_i$ is independent of the $\sigma$-algebra generated by $\{A_j:j\in I\setminus(N_D(i)\cup\{i\})\}$.

\begin{theorem}[{\cite{Spencer1977}, Theorem~1.5}]
Let $\mathcal A=\{A_1,\ldots,A_n\}$ be a finite family of events in a probability space. Suppose that $\Pr(A_i)\le p$ for every $i\in[n]$ and that the family $\mathcal A$ has a dependency graph of maximum degree at most $d$. If $ep(d+1)\le1$, then
$\Pr\Bigl(\bigcap_{i=1}^{n}\overline{A_i}\Bigr)>0$.
\end{theorem}

\begin{theorem}[{\cite{Spencer1977}, Theorem~1.1}]
Let $\mathcal A=\{A_i:i\in I\}$ be a finite family of events with dependency graph $D$. Suppose that there are real numbers $x_i\in(0,1)$ such that
$\Pr(A_i)\le x_i\prod_{j\in N_D(i)}(1-x_j)$ for every $i\in I$.
Then $\Pr\Bigl(\bigcap_{i\in I}\overline{A_i}\Bigr)>0$.
\end{theorem}

In the following proof, we use the Moser--Tardos resampling framework~\cite{MoserTardos2010} together with the variable-assignment approach developed by Harris~\cite{Harris2016}. Besides guaranteeing that the resampling algorithm avoids all bad events, this formulation bounds the probability of other events in the final assignment.

Let $\zeta_1,\ldots,\zeta_n$ be independent random variables, each with its prescribed set of possible values. A variable--value pair $(i,j)$ represents the condition $\zeta_i=j$. An \emph{atomic event} is a collection of variable--value pairs and therefore prescribes values for several variables. Thus an atomic event has the form
$E=\{\zeta_{i_1}=j_1,\ldots,\zeta_{i_t}=j_t\}$,
or equivalently
$E=\{(i_1,j_1),\ldots,(i_t,j_t)\}$.

Two pairs $(i,j)$ and $(i,j')$ \emph{conflict}, written $(i,j)\sim(i,j')$, if $j\ne j'$. If $z$ is a variable--value pair and $B$ is an atomic event, we write $z\sim B$ if $z$ conflicts with a pair contained in $B$. Let $\mathcal B$ be the family of atomic bad events under consideration. Let $E$ be an atomic event and let $Y\subseteq\mathcal B$. We say that $Y$ is \emph{orderable to $E$} if either $Y=\{E\}$ or there is an ordering $Y=\{B_1,\ldots,B_s\}$ such that, for every $i\in[s]$, some pair $z_i\in E$ satisfies
$z_i\sim B_i$ and 
$z_i\not\sim B_j$ for every $j<i$.
The empty set is also orderable to $E$. We say that $Y$ is \emph{assignable to $E$} if either $Y=\{E\}$ or there is an injective map $f:Y\to E$ such that every $B\in Y$ contains a variable--value pair conflicting with $f(B)$.

Harris~\cite{Harris2016} proved the following relationship between orderability and assignability.

\begin{lemma}[\cite{Harris2016}, Proposition~2.5]
\label{prop:harris-assignable}
If $Y$ is orderable to an atomic event $E$, then $Y$ is assignable to $E$.
\end{lemma}

The Moser--Tardos resampling algorithm starts from an independent random assignment. If no bad event occurs, it stops. Otherwise, it chooses an occurring bad event and independently resamples all variables involved in that event from their original distributions, leaving the remaining variables unchanged. The procedure is repeated until no bad event remains. 
Throughout the algorithm,
the event-selection rule depends only on the current and previous states. 
The variables of the
chosen event are independently resampled according to their original distributions.
We write $\PrOmega$ for probability under the initial independent assignment and $\PrMT$ for probability under the terminal distribution.

Given a weight function $\mu\colon\mathcal B\to[0,\infty)$, define
$\lambda(i,j):=
\sum_{j'\ne j}
\sum_{\substack{B'\in\mathcal B\\(i,j')\in B'}}
\mu(B')$
for a variable--value pair $(i,j)$, and define
$\psi(B):=\prod_{(i,j)\in B}\bigl(1+\lambda(i,j)\bigr)$
for an atomic event $B$. Harris~\cite{Harris2016} proved the following resampling criterion.

\begin{lemma}[\cite{Harris2016}, Proposition~2.7]
\label{prop:harris-general-criterion}
Let $\mathcal B$ be a family of atomic bad events and let $\mu\colon\mathcal B\to[0,\infty)$. Suppose that
\[
\mu(B)\ge\PrOmega(B)\bigl(\mu(B)+\psi(B)\bigr)
\qquad\text{for every }B\in\mathcal B.
\]
Then the Moser--Tardos algorithm terminates with probability one, and the expected number of resamplings of each $B\in\mathcal B$ is at most $\mu(B)$.
\end{lemma}

We also use Harris's estimate for the terminal distribution.

\begin{theorem}[\cite{Harris2016}, Theorem~4.3]
\label{thm:harris-terminal-general}
Suppose that $\mathcal B$ satisfies Harris's variable-assignment resampling criterion with weights $\mu$. If $E$ is an atomic event with $E\notin\mathcal B$, then
\[
\PrMT(E)
\le \PrOmega(E)
\sum_{\substack{Y\subseteq\mathcal B\\Y\text{ orderable to }E}}
\prod_{F\in Y}\mu(F).
\]
\end{theorem}

We now assume that all variables are independent Bernoulli random variables. Thus every atomic event $E$ can be written in the form
$E=\{\zeta_i=0\text{ for }i\in Z_E,\;
      \zeta_i=1\text{ for }i\in R_E\}$,
where $Z_E$ and $R_E$ are disjoint. So $Z_E$ is the set of indices of the variables that $E$ requires to be $0$, and $R_E$ is the set of indices of the variables that $E$ requires to be $1$. Let $\mathcal B$ be a finite family of distinct atomic bad events and assign each $E\in\mathcal B$ a positive weight $\mu(E)$. For each index $i$, define
\begin{equation}\label{eq:loads}
\Lambda_0(i):=\sum_{\substack{E\in\mathcal B\\i\in Z_E}}\mu(E),
\qquad
\Lambda_1(i):=\sum_{\substack{E\in\mathcal B\\i\in R_E}}\mu(E).
\end{equation}
For any atomic event $E$, define
\begin{equation}\label{eq:correction-factor}
\Phi(E):=
\prod_{i\in Z_E}\bigl(1+\Lambda_1(i)\bigr)
\prod_{i\in R_E}\bigl(1+\Lambda_0(i)\bigr).
\end{equation}

\begin{corollary}
\label{lem:harris-resampling}
Suppose that
\begin{equation}\label{eq:binary-resampling}
\PrOmega(E)\bigl(\mu(E)+\Phi(E)\bigr)\le\mu(E)
\qquad\text{for every }E\in\mathcal B.
\end{equation}
Then the Moser--Tardos algorithm terminates with probability one and its final assignment avoids every event in $\mathcal B$. Moreover, if $E\notin\mathcal B$ is atomic, then
\[
\PrMT(E)\le\PrOmega(E)\Phi(E).
\]
\end{corollary}

\begin{proof}
Choose an event $E\in\mathcal B$. Since the variables are binary, the product in \cref{prop:harris-general-criterion} satisfies
\begin{align*}
\psi(E)
&=\prod_{(i,j)\in E}
  \Biggl(1+\sum_{j'\ne j}
  \sum_{\substack{B'\in\mathcal B\\(i,j')\in B'}}\mu(B')\Biggr)\\
&=\prod_{i\in Z_E}\bigl(1+\Lambda_1(i)\bigr)
  \prod_{i\in R_E}\bigl(1+\Lambda_0(i)\bigr)\\
&=\Phi(E),
\end{align*}
where we have used \eqref{eq:loads} and \eqref{eq:correction-factor}. Thus \eqref{eq:binary-resampling} gives
$\PrOmega(E)\bigl(\mu(E)+\psi(E)\bigr)
=\PrOmega(E)\bigl(\mu(E)+\Phi(E)\bigr)
\le\mu(E)$.
By \Cref{prop:harris-general-criterion}, the Moser--Tardos algorithm terminates with probability one. At termination no bad event is present, so the terminal assignment avoids every event of $\mathcal B$.

Choose an atomic event $E\notin\mathcal B$. It remains to prove
$\PrMT(E)\le\PrOmega(E)\Phi(E)$.
By \cref{thm:harris-terminal-general},
$\PrMT(E)
\le\PrOmega(E)
\sum_{\substack{Y\subseteq\mathcal B\\Y\text{ orderable to }E}}
\prod_{F\in Y}\mu(F)$.
By \cref{prop:harris-assignable}, if $Y$ is orderable to $E$, then it is assignable to $E$. Hence
$\sum_{\substack{Y\subseteq\mathcal B\\Y\text{ orderable to }E}}
\prod_{F\in Y}\mu(F)
\le
\sum_{\substack{Y\subseteq\mathcal B\\Y\text{ assignable to }E}}
\prod_{F\in Y}\mu(F)$.
Now choose a subset $Y\subseteq\mathcal B$ that is assignable to $E$. Since $E\notin\mathcal B$, we have $Y\ne\{E\}$. So there is an injective map $f:Y\to E$ such that every $F\in Y$ conflicts with $f(F)$. Since $f$ is injective, each variable--value pair of $E$ can be assigned to at most one bad event. For each $(i,j)\in E$, we may select either zero or one bad event $F\ni(i,1-j)$. These account respectively for the terms
$1$ and $\sum_{\substack{F\in\mathcal B\\(i,1-j)\in F}}\mu(F)$.
Combining the terminal-distribution bound, assignability, and the preceding count, we obtain
\begin{align*}
\PrMT(E)
&\le\PrOmega(E)
  \sum_{\substack{Y\subseteq\mathcal B\\Y\text{ orderable to }E}}
  \prod_{F\in Y}\mu(F)\\
&\le\PrOmega(E)
  \sum_{\substack{Y\subseteq\mathcal B\\Y\text{ assignable to }E}}
  \prod_{F\in Y}\mu(F)\\
&\le\PrOmega(E)
  \prod_{i\in Z_E}\bigl(1+\Lambda_1(i)\bigr)
  \prod_{i\in R_E}\bigl(1+\Lambda_0(i)\bigr)\\
&=\PrOmega(E)\Phi(E),
\end{align*}
where the last equality is \eqref{eq:correction-factor}. This is the required bound.
\end{proof}

\section{Independent transversals and the main theorem}\label{sec:transversal}

\subsection{Structural reduction}

We use the following criterion of Thomassen~\cite{Thomassen1997,Thomassen1998}.

\begin{lemma}[{\cite{Thomassen1997,Thomassen1998}}]\label{lem:thomassen}
Let $G$ be a Hamiltonian graph with a fixed Hamilton cycle $C$. 
Color the edges of $C$ red and the remaining edges of $G$ green.
Suppose that
there exists a red-stable green-dominating set $S$ in $G$. 
Then $G$ contains a Hamilton cycle distinct from $C$.
\end{lemma}

We next give a general reduction for regular Hamiltonian graphs in terms of
the parameter $k=r-1$. Note that the following reduction applies to arbitrary $k$.

\begin{lemma}\label{lem:reduction}
Let $k\ge1$, and let $G$ be a $(k+1)$-regular graph with a
Hamilton cycle $C$. Let $H=G-E(C)$.
Color the edges of $C$ red and the remaining edges of $G$ green.
For each $v\in V(G)$, let $D_v=N_H[v]$, and let $\mathcal D$ be the family of distinct sets among the
$D_v$'s. Then
\begin{enumerate}[label=(\roman*)]
    \item every $D\in\mathcal D$ has size $k$;
    \item every vertex belongs to at most $k$ members of
    $\mathcal D$;
    \item every red-stable transversal of $\mathcal D$
    is a nonempty red-stable green-dominating set.
\end{enumerate}
Consequently, if $\mathcal D$ admits a red-stable
transversal, then $G$ contains a Hamilton cycle distinct from $C$.
\end{lemma}

\begin{proof}
We verify the three assertions separately. Since $C$ is a Hamilton cycle, every vertex of $G$ is incident with exactly two edges of $C$. The graph $G$ is $(k+1)$-regular, so deleting the edges of $C$ removes exactly two incident edges at every vertex. Hence $H=G-E(C)$ is $(k-1)$-regular. For each $v\in V(G)$, its closed neighborhood in $H$ therefore consists of $v$ itself together with its $k-1$ neighbors in $H$. Thus
$\lvert D_v\rvert=\lvert N_H[v]\rvert=1+d_H(v)=k$.
Every member of $\mathcal D$ is one of these sets $D_v$, so every member of $\mathcal D$ has size $k$. This proves (i).

For (ii), fix a vertex $y\in V(G)$. By definition, $y\in D_v=N_H[v]$ if and only if either $v=y$ or $v$ is adjacent to $y$ in $H$. Since $H$ is undirected, this is equivalent to $v\in N_H[y]$. The closed neighborhood $N_H[y]$ has size $k$ by (i). Hence exactly $k$ indices $v$ satisfy $y\in D_v$. The family $\mathcal D$ is obtained from the indexed family by deleting repeated sets, so identifying equal sets cannot increase the number of distinct members that contain $y$. Therefore $y$ belongs to at most $k$ members of $\mathcal D$, proving (ii).

For (iii), let $S$ be a red-stable transversal of $\mathcal D$. Since every indexed set $D_v$ is represented by a member of $\mathcal D$, the set $S$ meets $D_v$ for every $v\in V(G)$. Fix a vertex $v\notin S$. Because $S\cap D_v\ne\varnothing$ and $D_v=\{v\}\cup N_H(v)$, the intersection cannot occur at $v$ itself; hence $S\cap N_H(v)\ne\varnothing$. Thus some vertex of $S$ is adjacent to $v$ by an edge of $H$, that is, by a green edge. Since this holds for every $v\notin S$, the set $S$ is green-dominating. It is red-stable by hypothesis. Finally, $V(G)$ is nonempty and hence the indexed family $(D_v)$ is nonempty, so $\mathcal D$ is nonempty; a set meeting every member of a nonempty family cannot itself be empty. Therefore $S$ is a nonempty $C$-independent green-dominating set. This proves (iii).

By \cref{lem:thomassen}, 
$G$ has a Hamilton cycle distinct from the fixed cycle $C$. This proves the final assertion.
\end{proof}

\subsection{Probability space and estimates}
\label{sec:path}

We now prove \cref{thm:transversal}.
Let $k\ge17$, and let $C$ and $\mathcal D$ satisfy the hypotheses of
\cref{thm:transversal}. 
If $\mathcal D=\varnothing$, 
then the conclusion of \cref{thm:transversal} is obtained.
So assume $\mathcal D\neq\varnothing$.
Let $p=13/40$ and $q=27/40=1-p$.

Choose an arbitrary vertex $\rho\in V(C)$ and delete it from $C$. The
remaining graph is the path $P=C-\rho$.
Choose a bipartition $V(P)=A\mathbin{\dot\cup}B$ so that at most one endpoint of $P$ lies in $A$. This is possible by interchanging the names of the two parts if necessary.

Let $(\xi_b)_{b\in B}$ be mutually independent Bernoulli random
variables with $\Pr(\xi_b=1)=p$ and $\Pr(\xi_b=0)=q$. Define
$T:=\{b\in B:\xi_b=1\},
S:=T\cup\{a\in A:N_P(a)\cap T=\varnothing\}$.
In particular, $\rho\notin S$.

\begin{lemma}\label{lem:independent}
The random set $S$ is independent in the Hamilton cycle $C$.
\end{lemma}

\begin{proof}
Because $A\mathbin{\dot\cup}B$ is a bipartition of the path $P$, every edge of $P$ has one endpoint in $A$ and one endpoint in $B$. In particular, there are no edges inside $A$ and no edges inside $B$.

The set $T$ is a subset of $B$, so no two vertices of $T$ are adjacent in $P$. A vertex $a\in A$ is added to $S$ only when $N_P(a)\cap T=\varnothing$. Hence, if $b\in T$, then no neighbor of $b$ in $A$ belongs to $S$. Thus there is no edge of $P$ with one endpoint in $T$ and the other endpoint in $S\cap A$. Since there are also no edges inside either bipartition class, no edge of $P$ has both endpoints in $S$. Therefore $S$ is independent in $P$.

The only edges of the original cycle $C$ that are not present in $P=C-\rho$ are the two cycle edges incident with $\rho$. By construction $\rho\notin S$, so neither of those two edges can have both endpoints in $S$. Combining this with independence in $P$, we conclude that $S$ is independent in the whole cycle $C$.
\end{proof}

Order the vertices of $B$ as they appear from one end of $P$ to the other, and let $Q$ be the path whose vertex set is $B$ with consecutive vertices joined in this order. Thus $Q$ is an auxiliary path obtained by suppressing the vertices of $A$ from $P$; it is not a subgraph of $P$ in general. 

Choose a set $D\in\mathcal D$, and define
$Z_D:=D\cap B,
\mathcal F_D:=\{N_P(a):a\in D\cap A\}$.
If $a\in A$ is an internal vertex of $P$, then its two neighbors in $P$ are consecutive vertices of $B$ in the order defining $Q$, so $N_P(a)$ is a two-element set corresponding to an edge of $Q$. If $a$ is an endpoint of $P$, then $N_P(a)$ is a singleton containing the corresponding endpoint of $Q$. Our choice of the bipartition ensures that at most one endpoint of $P$ lies in $A$, so $\mathcal F_D$ contains at most one singleton. Also, two different vertices of $A$ have different neighborhoods in a path, so the members of $\mathcal F_D$ are distinct.

The vertices of $D\setminus\{\rho\}$ are partitioned into the vertices in $D\cap B$, counted by $\lvert Z_D\rvert$, and the vertices in $D\cap A$, each of which contributes one distinct member of $\mathcal F_D$. Therefore
$\lvert Z_D\rvert+\lvert\mathcal F_D\rvert
=\lvert D\setminus\{\rho\}\rvert$.
Let $\mathcal R_D$ be the family of minimal transversals of $\mathcal F_D$ that avoid $Z_D$. For $R\in\mathcal R_D$, define the atomic event
$E(D,R):=
\bigcap_{b\in Z_D}\{\xi_b=0\}
\cap\bigcap_{b\in R}\{\xi_b=1\}$.
On the one hand, the zero conditions ensure that the vertices of $D\cap B$ are not selected; 
on the other hand, the one conditions ensure that every vertex of $D\cap A$ has a selected neighbor in $T$. 
Therefore, no vertex of $D$ is added to $S$.

The next lemma shows that $S$ misses $D$ exactly 
when $E(D,R)$ occurs for some $R\in \mathcal R_D$.

\begin{lemma}\label{lem:encoding}
For every $D\in\mathcal D$,
\[
\{S\cap D=\varnothing\}=\bigcup_{R\in\mathcal R_D}E(D,R).
\]
\end{lemma}

\begin{proof}
We first show that
$\{S\cap D=\varnothing\}\subseteq
\bigcup_{R\in\mathcal R_D}E(D,R)$,
and then prove the reverse inclusion. First suppose that $S\cap D=\varnothing$. For a vertex $b\in B$, we have $b\in S$ if and only if $b\in T$. Hence $b\notin S$ if and only if $b\notin T$, which is equivalent to $\xi_b=0$. Since every vertex of $Z_D=D\cap B$ is outside $S$, we obtain $Z_D\cap T=\varnothing$.

Now consider a vertex $a\in D\cap A$. By definition, $a$ is added to $S$ exactly when none of its neighbors in $P$ lies in $T$. Since $a\notin S$, $N_P(a)\cap T\ne\varnothing$. Thus $T$ meets every member of $\mathcal F_D$. In other words, $T$ is a transversal of $\mathcal F_D$. Since the family is finite, $T$ has a minimal transversal $R\subseteq T$. Since $T$ is disjoint from $Z_D$, this minimal transversal also satisfies $R\cap Z_D=\varnothing$, so $R\in\mathcal R_D$. We have $\xi_b=0$ for every $b\in Z_D$ and $\xi_b=1$ for every $b\in R\subseteq T$. Therefore $E(D,R)$ occurs. 

Conversely, suppose that $E(D,R)$ occurs for some $R\in\mathcal R_D$. The zero conditions in $E(D,R)$ imply that $Z_D\cap T=\varnothing$. Hence no vertex of $D\cap B$ belongs to $S$. The one conditions imply $R\subseteq T$. 
Since $R$ is a transversal of $\mathcal F_D$, we have $N_P(a)\cap R\neq \varnothing$ for every $a\in D\cap A$. Therefore, $N_P(a)\cap T\ne\varnothing$. Thus no vertex of $D\cap A$ is added to $S$. Finally, if $\rho\in D$, then $\rho\notin S$ by construction. 
Therefore, we have shown that no vertex of $D$ can be added to $S$, and hence $S\cap D=\varnothing$. This proves the reverse inclusion.
\end{proof}

The preceding reformulation expresses the event that $S$ misses a member
$D\in\mathcal D$ in terms of minimal transversals of a family of path constraints.
For the following estimates, let $\alpha=29/40$ and $\beta=14/40$.
We next bound the total weight of these minimal transversals.

\begin{lemma}\label{lem:path-estimate}
Let $Q$ be a path, 
and let $\mathcal F$ be a family of subsets of $V(Q)$ such that each $F\in \mathcal F$ is the set of endvertices of an edge of $Q$, with the possible exception of one singleton consisting of an endvertex of $Q$. 
Let $\mathcal R(\mathcal F)$ be the family of minimal transversals of $\mathcal F$. For $b\in V(Q)$, define $d_{\mathcal F}(b):=\lvert\{F\in\mathcal F:b\in F\}\rvert$. 
Then
\[
\sum_{R\in\mathcal R(\mathcal F)}\beta^{\lvert R\rvert}
\le\alpha^{\lvert\mathcal F\rvert}.
\]
Moreover, for every $b\in V(Q)$,
\[
\sum_{\substack{R\in\mathcal R(\mathcal F)\\ b\in R}}
\beta^{\lvert R\rvert}
\le d_{\mathcal F}(b)\beta\alpha^{\lvert\mathcal F\rvert-1}.
\]
Both inequalities hold if the sums are restricted to any subfamily of $\mathcal R(\mathcal F)$, in particular for the transversals avoiding a prescribed set of vertices.
\end{lemma}

\begin{proof}

First, we prove the following facts, which will be used in the following proof.

\begin{claim}\label{clm:containing}
If $R$ is a minimal transversal of $\mathcal F$ and $x\in R$, then there exists $F_x\in\mathcal F$ such that $F_x\cap R=\{x\}$.
\end{claim}

\begin{proof}
Otherwise, every member of $\mathcal F$ containing $x$ would also contain another vertex of $R$. 
Hence $R\setminus\{x\}$ would still meet every member of $\mathcal F$, contrary to the minimality of $R$.  
\end{proof}

\begin{claim}\label{clm:delete}
For a vertex $b$, let $\mathcal F_b:=\{F\in\mathcal F:b\notin F\}$. Suppose that $R\in\mathcal R(\mathcal F)$ and $b\in R$. Then $R\setminus\{b\}$ is a minimal transversal of $\mathcal F_b$.
\end{claim}

\begin{proof}
If $F\in\mathcal F_b$, then $b\notin F$. 
Since $R$ meets $F$, $F\cap (R\setminus\{b\})\neq\varnothing$.
So $R\setminus\{b\}$ meets every member of $\mathcal F_b$.

It remains to prove minimality. 
Let $x\in R\setminus\{b\}$. 
By \cref{clm:containing}, there exists $F_x\in\mathcal F$ such that $F_x\cap R=\{x\}$. Since $b\in R$ and $b\ne x$, we have $b\notin F_x$, so $F_x\in\mathcal F_b$ and $F_x\cap(R\setminus\{b\})=\{x\}$. Thus removing $x$ would leave $F_x$ uncovered. Hence $R\setminus\{b\}$ is a minimal transversal of $\mathcal F_b$. 
\end{proof}

We now prove the first inequality. 
First suppose that $\mathcal F$ contains no singleton. 
Let $H$ be the spanning subgraph of $Q$ consisting of the edges whose endvertex sets belong to $\mathcal F$. Thus each component of $H$ is a path or an isolated vertex.
By the definition of $H$, a transversal of $\mathcal F$ is exactly a vertex cover of this graph. 
Moreover, a minimal vertex cover contains no isolated vertex, and its intersection with each nontrivial component is a minimal vertex cover of that component. 
Conversely, the union of one minimal vertex cover from each nontrivial component is a minimal vertex cover of the whole graph. 
Since the nontrivial components of $H$ are vertex-disjoint, every minimal transversal $R$ is uniquely the union of minimal vertex covers $R_j$ of these components. 
Moreover, $\beta^{\lvert R\rvert}=\prod_j\beta^{\lvert R_j\rvert}$. 
Hence the weighted sum is the product of the corresponding sums for the nontrivial path components, so it is enough to consider one path.

For $\ell\ge0$, let
$g_\ell=\sum_R\beta^{\lvert R\rvert}$,
where the sum is taken over all minimal vertex covers $R$ of a path with $\ell$ edges.  
For $\ell=0$, the empty set is the unique minimal cover, so $g_0=\beta^{\lvert\varnothing\rvert}=1$. For $\ell=1$, either endpoint alone is a minimal cover, so $g_1=2\beta$. For $\ell=2$, the minimal covers are the middle vertex and the pair of endvertices, so $g_2=\beta+\beta^2$.

Let $\ell\ge3$, and let $u,v,w$ be the first three vertices of the path, with $u$ an endvertex. We divide the minimal vertex covers according to whether they contain $u$.

Suppose first that $u\notin R$. Since $uv$ must be covered, $v\in R$. After deleting $u$ and $v$ from the path, the set $R\setminus\{v\}$ covers the remaining path, which has $\ell-2$ edges. 
We now show that $R\setminus \{v\}$ is a minimal vertex cover of the remaining path. Indeed, if $x\in R\setminus\{v\}$, by \cref{clm:containing}, we can choose an edge $F_x$ with $F_x\cap R=\{x\}$. 
Since $v\in R$ and $v\neq x$, the edge $F_x$ is not incident to $v$. 
Since $u$ is an endvertex of the path, $F_x$ remains in the remaining path and still meets $R\setminus\{v\}$ only in $x$. 
So $F_x\cap (R\setminus \{v\})=\{x\}$. If we delete $x$ from $R\setminus \{v\}$, then $F_x$ is uncovered.
Thus, $R\setminus \{v\}$ is a minimal vertex cover of the remaining path.
Conversely, if $R'$ is a minimal vertex cover of the shorter path, then $R'\cup\{v\}$ is a minimal vertex cover of the original path. 
If we delete $v$ from $R'\cup\{v\}$, then $uv$ is uncovered.
Since the remaining path is unchanged on those edges, no vertex of $R'$ can be deleted.
Therefore the total weight of the minimal vertex covers with $u\notin R$ is $\beta g_{\ell-2}$.

Suppose next that $u\in R$. Then $v\notin R$; otherwise removing $u$ would leave $uv$ covered. Since $v\notin R$, the edge $vw$ forces $w\in R$. After deleting $u,v,w$, the set $R\setminus\{u,w\}$ covers the remaining path, which has $\ell-3$ edges. 
Now we show that $R\setminus\{u,w\}$ is a minimal vertex cover of the remaining path. To see this, let $x\in R\setminus\{u,w\}$. By \cref{clm:containing}, we can choose an edge $F_x$ with $F_x\cap R=\{x\}$.
Since $u,w\in R$ and $x\notin\{u,w\}$, the edge $F_x$ is incident to neither $u$ nor $w$. 
The only edges incident to $v$ are $uv$ and $vw$,
so $F_x$ lies in the remaining path and prevents $x$ from being removed. 
So $R\setminus\{u,w\}$ is a minimal vertex cover of the remaining path.
Conversely, if $R'$ is a minimal vertex cover of the remaining path, then $R'\cup\{u,w\}$ is a minimal vertex cover of the original path,
since removing $u$ uncovers $uv$, removing $w$ uncovers $vw$, and every vertex of $R'$ remains necessary on the remaining path. Thus the covers with $u\in R$ contribute $\beta^2g_{\ell-3}$. Therefore
$g_\ell=\beta g_{\ell-2}+\beta^2g_{\ell-3}$.

We prove $g_\ell\le\alpha^\ell$ by induction on $\ell$. For the base cases, we have $g_0=1$. 
Since
$\alpha-2\beta=\frac{1}{40},
\alpha^2-\beta-\beta^2=\frac{17}{320}$,
we have $g_1=2\beta\le\alpha$ and $g_2=\beta+\beta^2\le\alpha^2$.
Let $\ell\ge3$ and assume the bound holds for all smaller indices.
The recurrence and the induction hypothesis give
\begin{align*}
g_\ell
&\le\beta\alpha^{\ell-2}+\beta^2\alpha^{\ell-3}=\alpha^{\ell-3}(\alpha\beta+\beta^2)\le\alpha^\ell,
\end{align*}
where the last inequality follows from
$\alpha^3-\alpha\beta-\beta^2=309/64000>0$.
Thus $g_\ell\le\alpha^\ell$ for every $\ell\ge0$. 
Let $H_1,\ldots,H_t$ be the nontrivial components of $H$. For each $j$, let $R_j$ range over the minimal vertex covers of $H_j$.
Let $\ell_j=\lvert E(H_j)\rvert$ for $j=1,2,\ldots,t$.
So $\beta^{\lvert R\rvert}=\prod_{j=1}^t\beta^{\lvert R_j\rvert}$.
Therefore,
\begin{align*}
\sum_{R\in\mathcal R(\mathcal F)}\beta^{\lvert R\rvert}
&=\sum_{R_1}\cdots\sum_{R_t}
  \prod_{j=1}^t\beta^{\lvert R_j\rvert}\\
&=\prod_{j=1}^t\biggl(\sum_{R_j}\beta^{\lvert R_j\rvert}\biggr)
 =\prod_{j=1}^tg_{\ell_j}\\
&\le\prod_{j=1}^t\alpha^{\ell_j}
 =\alpha^{\lvert\mathcal F\rvert}.
\end{align*}

Now suppose that $\mathcal F$ also contains the allowed singleton $\{u\}$, where $u$ is an endvertex of $Q$. Every transversal of $\mathcal F$ contains $u$. For each $R\in\mathcal R(\mathcal F)$, by \cref{clm:delete}, $R\setminus\{u\}$ is a minimal transversal of $\mathcal F_u$. 
If $R_1$ and $R_2$ are different sets, then $R_1\setminus\{u\}$ and $R_2\setminus \{u\}$ are different sets. 
Since the family $\mathcal F_u$ contains no singleton, it has $\lvert\mathcal F\rvert-d_{\mathcal F}(u)$ members. 
Since $u$ is an endvertex of $Q$, it belongs to the singleton and to at most one edge, so $d_{\mathcal F}(u)\in\{1,2\}$. Hence
\begin{align*}
\sum_{R\in\mathcal R(\mathcal F)}\beta^{\lvert R\rvert}
&\le\beta\sum_{R'\in\mathcal R(\mathcal F_u)}\beta^{\lvert R'\rvert}\le\beta\alpha^{\lvert\mathcal F\rvert-d_{\mathcal F}(u)}.
\end{align*}
If $d_{\mathcal F}(u)=1$, this is at most $\alpha^{\lvert\mathcal F\rvert}$ because $\beta\le\alpha$. If $d_{\mathcal F}(u)=2$, it is at most $\alpha^{\lvert\mathcal F\rvert}$ because $\beta\le\alpha^2$. This proves the first inequality in all cases.

We now prove the second bound. Choose a vertex $b\in V(Q)$, and let $d=d_{\mathcal F}(b)$. 
By the definition of $\mathcal F$, we have $d\in\{0,1,2\}$. 
If $d=0$, then no minimal transversal of $\mathcal F$ contains $b$, 
so
$\sum_{\substack{R\in\mathcal R(\mathcal F)\\ b\in R}}
\beta^{\lvert R\rvert}
=d_{\mathcal F}(b)\beta\alpha^{\lvert\mathcal F\rvert-1}=0$.
Now we assume that $d>0$.
By \cref{clm:delete}, for every $R\in\mathcal R(\mathcal F)$ containing $b$, the set $R\setminus\{b\}$ is a minimal transversal of $\mathcal F_b$.
The map $R\mapsto R\setminus\{b\}$ is injective. The first bound for $\mathcal F_b$ yields
\begin{align*}
\sum_{\substack{R\in\mathcal R(\mathcal F)\\ b\in R}}
\beta^{\lvert R\rvert}
&\le\beta\sum_{R'\in\mathcal R(\mathcal F_b)}
\beta^{\lvert R'\rvert}\le\beta\alpha^{\lvert\mathcal F\rvert-d}\le d\beta\alpha^{\lvert\mathcal F\rvert-1}.
\end{align*}
For $d=1$, the last inequality is an equality. For $d=2$, since $1\le2\alpha$, the last inequality holds.

Finally, let $X\subseteq V(Q)$, and let
$\mathcal R_X(\mathcal F)
=\{R\in\mathcal R(\mathcal F):R\cap X=\varnothing\}$.
Since $\mathcal R_X(\mathcal F)\subseteq\mathcal R(\mathcal F)$ and all summands are nonnegative, 
both estimates hold with $\mathcal R(\mathcal F)$ replaced by $\mathcal R_X(\mathcal F)$.
\end{proof}

Let
$\mathcal D_0:=\{D\in\mathcal D:\rho\notin D\},
\mathcal D_1:=\{D\in\mathcal D:\rho\in D\}$.
Since every vertex belongs to at most $k$ members of $\mathcal D$, we have
$\lvert\mathcal D_1\rvert\le k$. Let $\mathcal B$ be the family of distinct
atomic events $E(D,R)$ with $D\in\mathcal D_0$ and $R\in\mathcal R_D$, and
assign the weight
$\mu\bigl(E(D,R)\bigr):=\alpha^{\lvert Z_D\rvert}\beta^{\lvert R\rvert}$.
If two pairs $(D_1,R_1)$ and $(D_2,R_2)$ define the same atomic event, then
they require the same variables to be $0$ and the same variables to be $1$.
Hence they have the same $\mu$-value. Recall that for $b\in B$, $\Lambda_0(b)$
and $\Lambda_1(b)$ denote the total weights of the bad events requiring
$\xi_b=0$ and $\xi_b=1$, respectively. We next obtain uniform bounds on
$\Lambda_0(b)$ and $\Lambda_1(b)$, which will be used to verify the hypothesis
of \cref{lem:harris-resampling}.

\begin{lemma}\label{lem:loads}
For every $b\in B$, we have
\[
\Lambda_0(b)\le k\alpha^k<\frac{2}{27},
\qquad
\Lambda_1(b)\le2k\beta\alpha^{k-1}<\frac{1}{14}.
\]
Moreover, the sequences $k\alpha^k$, $2k\beta\alpha^{k-1}$, and
$k\alpha^{k-1}$ are decreasing for $k\ge17$.
\end{lemma}

\begin{proof}
Choose a vertex $b\in B$. We first consider $\Lambda_0(b)$. If $E(D,R)$ is
included in the sum defining $\Lambda_0(b)$, then $b\in Z_D$. Since
$Z_D=D\cap B$, we have $b\in D$. By the hypothesis of \cref{thm:transversal},
the vertex $b$ belongs to at most $k$ members of $\mathcal D$, and therefore
to at most $k$ members of the subfamily $\mathcal D_0$.

Fix such a set $D$. By \cref{lem:path-estimate},
$\sum_{R\in\mathcal R_D}\beta^{\lvert R\rvert}
\le\alpha^{\lvert\mathcal F_D\rvert}$.
Since $D\in\mathcal D_0$ does not contain $\rho$, we also have
$\lvert Z_D\rvert+\lvert\mathcal F_D\rvert=\lvert D\rvert=k$. Thus the total
weight of the events arising from $D$ satisfies
\begin{align*}
\sum_{R\in\mathcal R_D}\mu\bigl(E(D,R)\bigr)
&=\alpha^{\lvert Z_D\rvert}
  \sum_{R\in\mathcal R_D}\beta^{\lvert R\rvert}\le\alpha^{\lvert Z_D\rvert}\alpha^{\lvert\mathcal F_D\rvert}
 =\alpha^k.
\end{align*}
Since $b$ belongs to at most $k$ members of $\mathcal D_0$, it follows that
\begin{align*}
\Lambda_0(b)
&\le\sum_{\substack{D\in\mathcal D_0\\b\in D}}
     \sum_{R\in\mathcal R_D}\mu\bigl(E(D,R)\bigr)\le\sum_{\substack{D\in\mathcal D_0\\b\in D}}\alpha^k
 \le k\alpha^k.
\end{align*}

We next consider $\Lambda_1(b)$. If $E(D,R)$ is included in the sum defining
$\Lambda_1(b)$, then $b\in R$. For a fixed $D\in\mathcal D_0$ with
$d_{\mathcal F_D}(b)>0$, \cref{lem:path-estimate} gives
$\sum_{\substack{R\in\mathcal R_D\\b\in R}}\beta^{\lvert R\rvert}
\le d_{\mathcal F_D}(b)\beta\alpha^{\lvert\mathcal F_D\rvert-1}$.
Hence
\begin{align*}
\sum_{\substack{R\in\mathcal R_D\\b\in R}}\mu\bigl(E(D,R)\bigr)
&=\alpha^{\lvert Z_D\rvert}
  \sum_{\substack{R\in\mathcal R_D\\b\in R}}\beta^{\lvert R\rvert}\\
&\le d_{\mathcal F_D}(b)\beta
     \alpha^{\lvert Z_D\rvert+\lvert\mathcal F_D\rvert-1}\\
&=d_{\mathcal F_D}(b)\beta\alpha^{k-1}.
\end{align*}

It remains to estimate $\sum_{D\in\mathcal D_0}d_{\mathcal F_D}(b)$.
By definition,
$d_{\mathcal F_D}(b)=\sum_{a\in N_P(b)}1_{\{a\in D\}}$.
Therefore
\begin{align*}
\sum_{D\in\mathcal D_0}d_{\mathcal F_D}(b)
&=\sum_{D\in\mathcal D_0}\sum_{a\in N_P(b)}1_{\{a\in D\}}\\
&=\sum_{a\in N_P(b)}
  \bigl\lvert\{D\in\mathcal D_0:a\in D\}\bigr\rvert\\
&\le k\lvert N_P(b)\rvert\le2k.
\end{align*}
Consequently,
\begin{align*}
\Lambda_1(b)
&\le\sum_{D\in\mathcal D_0}
     \sum_{\substack{R\in\mathcal R_D\\b\in R}}\mu\bigl(E(D,R)\bigr)\\
&\le\beta\alpha^{k-1}\sum_{D\in\mathcal D_0}d_{\mathcal F_D}(b)\\
&\le2k\beta\alpha^{k-1}.
\end{align*}

For each of the sequences $k\alpha^k$, $2k\beta\alpha^{k-1}$, and
$k\alpha^{k-1}$, the ratio of two consecutive terms satisfies
$\frac{k+1}{k}\alpha
\le\frac{18}{17}\cdot\frac{29}{40}
=\frac{261}{340}<1 (k\ge17)$.
Hence all three sequences are decreasing for $k\ge17$.

Since $\alpha^2=841/1600<10/19$, we have
$\alpha^8<\left(\frac{10}{19}\right)^4<\frac{1}{13},
\alpha^{16}<\frac{1}{169}$.
It follows that
$17\alpha^{17}
<17\cdot\frac{29}{40}\cdot\frac{1}{169}<\frac{2}{27},
34\beta\alpha^{16}
<34\cdot\frac{14}{40}\cdot\frac{1}{169}<\frac{1}{14}$.
Thus, for every $k\ge17$, we have $k\alpha^k<2/27$ and
$2k\beta\alpha^{k-1}<1/14$.
\end{proof}

\begin{lemma}\label{lem:resampling-check}
The weights $\mu\bigl(E(D,R)\bigr)=\alpha^{\lvert Z_D\rvert}
\beta^{\lvert R\rvert}$ satisfy the hypothesis of
\cref{lem:harris-resampling} for every event in $\mathcal B$.
\end{lemma}

\begin{proof}
Choose an event $E=E(D,R)\in\mathcal B$, and write $z=\lvert Z_D\rvert$ and
$r=\lvert R\rvert$. We must verify
$\PrOmega(E)\bigl(\mu(E)+\Phi(E)\bigr)\le\mu(E)$.
Since $\mu\bigl(E(D,R)\bigr)=\alpha^{\lvert Z_D\rvert}
\beta^{\lvert R\rvert}>0$, this is equivalent to
$\PrOmega(E)+\frac{\PrOmega(E)\Phi(E)}{\mu(E)}\le1$.
We bound the two terms separately.

Since $R$ is a transversal of $\mathcal F_D$, every member of
$\mathcal F_D$ contains a vertex of $R$. Moreover, each $b\in R$ belongs to
at most two members of $\mathcal F_D$. Hence
$\lvert\mathcal F_D\rvert
\le\sum_{b\in R}d_{\mathcal F_D}(b)
\le2\lvert R\rvert$.
Since $D\in\mathcal D_0$, we have
$k=\lvert Z_D\rvert+\lvert\mathcal F_D\rvert
\le\lvert Z_D\rvert+2\lvert R\rvert$.
In particular,
$\lvert Z_D\rvert+\lvert R\rvert\ge\lceil k/2\rceil\ge9$.

By independence, $\PrOmega(E)=q^{\lvert Z_D\rvert}p^{\lvert R\rvert}$.
Since $p=13/40<(27/40)^2=q^2$ and
$\lvert Z_D\rvert+2\lvert R\rvert\ge k$, we obtain
$\PrOmega(E)
=q^{\lvert Z_D\rvert}p^{\lvert R\rvert}
\le q^{\lvert Z_D\rvert+2\lvert R\rvert}
 \le q^k\le q^{17}$.
Moreover, $q^3=(27/40)^3<5/16$, so
$q^{12}<\left(\frac{5}{16}\right)^4<\frac{1}{64}$.
Hence $\PrOmega(E)\le q^{17}<q^{12}<1/64$.

By \cref{lem:loads}, we have $\Lambda_1(b)<1/14$ for $b\in Z_D$, and
therefore
$$\frac{q\bigl(1+\Lambda_1(b)\bigr)}{\alpha}
<\frac{(27/40)(15/14)}{29/40}
=\frac{405}{406}.$$
Similarly, $\Lambda_0(b)<2/27$ for $b\in R$, and hence
$$\frac{p\bigl(1+\Lambda_0(b)\bigr)}{\beta}
<\frac{(13/40)(29/27)}{14/40}
=\frac{377}{378}<\frac{405}{406}.$$
Using the definitions of $E$, $\Phi$, and $\mu$, followed by these bounds
and $z+r\ge9$, we obtain
\begin{align*}
\frac{\PrOmega(E)\Phi(E)}{\mu(E)}
&=\prod_{b\in Z_D}\frac{q\bigl(1+\Lambda_1(b)\bigr)}{\alpha}
  \prod_{b\in R}\frac{p\bigl(1+\Lambda_0(b)\bigr)}{\beta}\\
&<\left(\frac{405}{406}\right)^{\lvert Z_D\rvert+\lvert R\rvert}\\
&\le\left(\frac{405}{406}\right)^9\\
&<\frac{405}{414}=\frac{45}{46}<\frac{63}{64}.
\end{align*}
Here by Bernoulli's inequality, we have 
$(406/405)^9 = (1 + 1/405)^9 > 1 + 9/405 = 414/405$.
So $(405/406)^9 < 405/414$.
Consequently,
$\PrOmega(E)+\frac{\PrOmega(E)\Phi(E)}{\mu(E)}
<\frac{1}{64}+\frac{63}{64}=1$.
Since $\mu(E)>0$, we have
$\PrOmega(E)\bigl(\mu(E)+\Phi(E)\bigr)\le\mu(E)$,
as required by \cref{lem:harris-resampling}.
\end{proof}

\begin{proof}[Proof of \cref{thm:transversal}]
If $\mathcal D=\varnothing$, the empty set meets every member of
$\mathcal D$ trivially and is independent in $C$, so it is a red-stable
transversal. We therefore assume $\mathcal D\ne\varnothing$.

Let $S$ be the random set defined in the preceding construction. By
\cref{lem:independent}, $S$ is red-stable for every realization of the
Bernoulli variables. Hence it remains to find a terminal assignment such
that $S\cap D\ne\varnothing$ for every $D\in\mathcal D$. Equivalently, it
suffices to ensure the event
$\bigcap_{D\in\mathcal D}\{S\cap D\ne\varnothing\}$.

Partition $\mathcal D$ as $\mathcal D_0\mathbin{\dot\cup}\mathcal D_1$, where
$\mathcal D_0$ consists of the sets avoiding $\rho$ and $\mathcal D_1$
consists of the sets containing $\rho$. We first enforce all constraints
from $\mathcal D_0$ exactly. By \cref{lem:encoding}, for a fixed
$D\in\mathcal D_0$ the event $S\cap D=\varnothing$ is the union of the atomic
events $E(D,R)$ with $R\in\mathcal R_D$. Every distinct atomic event arising
in this way belongs to the bad-event family $\mathcal B$. By
\cref{lem:resampling-check}, the weights assigned to $\mathcal B$ satisfy
the hypothesis of \cref{lem:harris-resampling}. Hence the Moser--Tardos
algorithm terminates with probability one and the terminal assignment
avoids every event in $\mathcal B$. If some $D\in\mathcal D_0$ were missed
by the terminal set, then \cref{lem:encoding} would imply that at least one
event $E(D,R)\in\mathcal B$ occurs, contradicting termination. Thus every
terminal assignment produced by the algorithm meets every member of
$\mathcal D_0$.

Fix $D\in\mathcal D_1$. By \cref{lem:encoding}, if $S\cap D=\varnothing$,
then at least one atomic event $E(D,R)$ with $R\in\mathcal R_D$ occurs. We
estimate the terminal probability of each such atomic event. It can happen
that the same atomic specification also arises from some pair $(D',R')$
with $D'\in\mathcal D_0$; in that case the distinct event $E(D,R)$ already
belongs to $\mathcal B$ and is avoided by every terminal assignment, so
its terminal probability is $0$. Otherwise $E(D,R)\notin\mathcal B$, and
the second conclusion of \cref{lem:harris-resampling} applies.

For such an event, its initial probability is
$q^{\lvert Z_D\rvert}p^{\lvert R\rvert}$, while its terminal correction
factor is
$\prod_{b\in Z_D}\bigl(1+\Lambda_1(b)\bigr)
\prod_{b\in R}\bigl(1+\Lambda_0(b)\bigr)$.
The load bounds in \cref{lem:loads} imply
$q\bigl(1+\Lambda_1(b)\bigr)<\alpha,
p\bigl(1+\Lambda_0(b)\bigr)<\beta$
for every coordinate. Thus, in both cases, we obtain
$\PrMT\bigl(E(D,R)\bigr)
\le\alpha^{\lvert Z_D\rvert}\beta^{\lvert R\rvert}$.

Applying the union bound and \cref{lem:encoding}, followed by the preceding
event bound and \cref{lem:path-estimate}, gives
\begin{align*}
\PrMT\bigl(S\cap D=\varnothing\bigr)
&\le\sum_{R\in\mathcal R_D}\PrMT\bigl(E(D,R)\bigr)\le\alpha^{\lvert Z_D\rvert}
     \sum_{R\in\mathcal R_D}\beta^{\lvert R\rvert}\le\alpha^{\lvert Z_D\rvert+\lvert\mathcal F_D\rvert}=\alpha^{k-1}.
\end{align*}
The last equality uses $\rho\in D$ and
$\lvert Z_D\rvert+\lvert\mathcal F_D\rvert
=\lvert D\setminus\{\rho\}\rvert=k-1$.
Thus every individual set in $\mathcal D_1$ is missed with probability at
most $\alpha^{k-1}$.

There are at most $k$ such sets. Indeed, $\mathcal D_1$ is precisely the
subfamily of members containing $\rho$, and by hypothesis every vertex,
in particular $\rho$, belongs to at most $k$ members of $\mathcal D$. A
second union bound, together with $\lvert\mathcal D_1\rvert\le k$, the
monotonicity in \cref{lem:loads}, and $\alpha^{16}<1/169$, gives
\begin{align*}
\PrMT\bigl(\exists D\in\mathcal D_1:S\cap D=\varnothing\bigr)
&\le\lvert\mathcal D_1\rvert\alpha^{k-1}\le k\alpha^{k-1}\le17\alpha^{16}<\frac{17}{169}<1.
\end{align*}
Because this failure probability is strictly less than one, the
complementary event has positive terminal probability. Hence there exists
at least one terminal assignment that meets every member of $\mathcal D_1$.
Every terminal assignment already meets every member of $\mathcal D_0$,
so this particular assignment meets every member of $\mathcal D$.

Let $S$ be the set corresponding to such an assignment. By
\cref{lem:independent}, $S$ is independent in $C$, and by the preceding
argument it intersects every member of $\mathcal D$. Therefore $S$ is a
red-stable transversal of $\mathcal D$, as required.
\end{proof}

\subsection{Proof of the main theorem}
\label{sec:calculations}

\begin{proof}[Proof of \cref{thm:main}]
Let $G$ be an $r$-regular graph with $r\geq 17$, and fix one
Hamilton cycle $C$ of $G$. We shall find a Hamilton cycle different from
this fixed cycle.

If $r=17$, then every vertex of $G$ has odd degree. Thomason~\cite{Thomason1978}
proved that a Hamiltonian graph in which every vertex has odd degree cannot
have a unique Hamilton cycle. Hence the desired second Hamilton cycle
exists when $r=17$.

We may therefore assume $r\ge18$. Let $k=r-1$, so $k\ge17$. Color the edges
of the fixed Hamilton cycle $C$ red and all remaining edges green, as above,
and let $H=G-E(C)$ be the spanning graph of green edges. Since every vertex
loses exactly its two incident cycle edges, $H$ is $(r-2)=(k-1)$-regular.
For each vertex $v$, let $D_v=N_H[v]$, and let $\mathcal D$ be the family of
distinct sets among these closed green neighborhoods.

By \cref{lem:reduction}, each member of $\mathcal D$ has exactly $k$ vertices
and every vertex of $C$ belongs to at most $k$ members of $\mathcal D$. Thus
the pair $(C,\mathcal D)$ satisfies every hypothesis of
\cref{thm:transversal}. Since $k\ge17$, that theorem supplies a red-stable
transversal $S$ of $\mathcal D$.

We now translate this transversal back to the graph problem. The set $S$
is independent on the Hamilton cycle $C$ by definition. Moreover, because
it meets every closed green neighborhood $D_v$, if a vertex $v$ is not
itself in $S$, then the intersection $S\cap D_v$ must contain a green
neighbor of $v$. Hence every vertex outside $S$ has a green neighbor in
$S$. In other words, $S$ is green-dominating. The family $\mathcal D$ is
nonempty, so its transversal $S$ is nonempty. Therefore $S$ is a nonempty
red-stable green-dominating set.

Therefore, by \cref{lem:thomassen}, $G$ has a Hamilton cycle distinct from
$C$. This proves the theorem for every $r\ge18$, and together with the
$r=17$ case proves it for every $r>16$.
\end{proof}

\section{Constructing a red-stable transversal}
\label{sec:algorithm}

We give the algorithms after completing the existence proofs. Let $n=\lvert V(C)\rvert$ and $m=\lvert\mathcal D\rvert$. Algorithms~\ref{alg:transversal} and~\ref{alg:minimal-witness} implement the resampling construction. 

\subsection{The algorithms}

Fix $\rho\in V(C)$, let $P=C-\rho$, and let $A\mathbin{\dot\cup}B$ be the bipartition with at most one endvertex in $A$. Let $\mathcal D_0$ and $\mathcal D_1$ consist of the sets avoiding and containing $\rho$, respectively. For $D\in\mathcal D$, let $Z_D=D\cap B$ and $\mathcal F_D=\{N_P(a):a\in D\cap A\}$. These are defined in \Cref{sec:path}. Algorithm~\ref{alg:transversal} uses Algorithm~\ref{alg:minimal-witness} to find a true atomic event whenever a set is missed. 
Fix all ordering rules in advance, including the rule for 
choosing among several available candidates. 
At each restart, reset the event-selection history
and reconstruct all data for the current attempt from 
the new initial assignment, according to these fixed rules.

\begin{algorithm}[H]
\caption{An algorithm for a red-stable transversal}
\label{alg:transversal}
\small
\textbf{Input:} $C,\mathcal D,k$ satisfying \cref{thm:transversal}.\par
\textbf{Output:} $S\subseteq V(C)$ independent in $C$ and meeting every $D\in\mathcal D$.\par\smallskip
\begin{tabularx}{\linewidth}{@{}rL@{}}
\algline{1}{If $\mathcal D=\varnothing$, \textbf{return} $\varnothing$.}
\algline{2}{Choose $\rho$ and construct $P,A,B,\mathcal D_0,\mathcal D_1,Z_D,\mathcal F_D$ as above.}
\algline{3}{Build incidence lists $\mathcal I(v)=\{D\in\mathcal D:v\in D\}$.}
\algline{4}{\textbf{repeat}}
\algline{5}{\quad Draw fresh independent $\xi_b\sim\mathrm{Bernoulli}(13/40)$ for all $b\in B$.}
\algline{6}{\quad Set $T=\{b\in B:\xi_b=1\}$ and $S=T\cup\{a\in A:N_P(a)\cap T=\varnothing\}$.}
\algline{7}{\quad Compute $h(D)=|S\cap D|$ and $\mathcal M=\{D\in\mathcal D_0:h(D)=0\}$.}
\algline{8}{\quad \textbf{while} $\mathcal M\ne\varnothing$ \textbf{do}}
\algline{9}{\qquad Choose $D\in\mathcal M$ using the fixed rule and only the history of the current attempt.}
\algline{10}{\qquad $R\gets\textsc{MinimalWitness}(\mathcal F_D,T)$. }
\algline{11}{\qquad Set $U=Z_D\cup R$ and independently redraw $\xi_b$ for every $b\in U$.}
\algline{12}{\qquad Update $T$, and membership in $S$ only on $U\cup N_P(U)$.}
\algline{13}{\qquad Using $\mathcal I(v)$, update affected counts $h(D')$ and maintain $\mathcal M$.}
\algline{14}{\quad \textbf{end while}}
\algline{15}{\quad \textbf{if} $h(D)>0$ for every $D\in\mathcal D_1$ \textbf{then return} $S$.}
\algline{16}{\textbf{end repeat} }
\end{tabularx}
\end{algorithm}

\begin{algorithm}[H]
\caption{\textsc{MinimalWitness}$(\mathcal F,T)$}
\label{alg:minimal-witness}
\small
\textbf{Input:} $T$ hits every $F\in\mathcal F$; $\mathcal F$ consists of edges of a path, with at most one singleton at an end-vertex.\par
\textbf{Output:} A minimal transversal $R\subseteq T$ of $\mathcal F$.\par\smallskip
\begin{tabularx}{\linewidth}{@{}rL@{}}
\algline{1}{Build a list $L=T\cap\bigcup_{F\in\mathcal F}F$ without repeated vertices; set $R\gets L$.}
\algline{2}{For every $F\in\mathcal F$, set $c(F)=|R\cap F|$.}
\algline{3}{Build $J(b)=\{F\in\mathcal F:b\in F\}$ for each $b\in L$.}
\algline{4}{\textbf{for each} $b$ in $L$, visited once, \textbf{do}}
\algline{5}{\quad \textbf{if} $c(F)\ge2$ for every $F\in J(b)$ \textbf{then}}
\algline{6}{\qquad Delete $b$ from $R$ and decrease $c(F)$ by $1$ for every $F\in J(b)$.}
\algline{7}{\quad \textbf{end if}}
\algline{8}{\textbf{end for}}
\algline{9}{\textbf{return} $R$.}
\end{tabularx}
\end{algorithm}

\subsection{Correctness and running time}

In Algorithm~\ref{alg:transversal}, maintain $\mathcal M$ using an array and a position map, allowing constant-time insertion, deletion, and selection. 
A set is removed from $\mathcal M$ only when its count becomes positive. At each iteration, 
select the first set in $\mathcal M$.

Algorithm~\ref{alg:minimal-witness} maintains $R\cap F\neq\varnothing$ for every $F\in \mathcal F$, since a vertex is deleted only when every constraint containing it remains covered.
If a vertex $b$ is not deleted, then some $F_b\in\mathcal F$ satisfies $R\cap F_b=\{b\}$ when $b$ is examined. 
Since the algorithm only deletes vertices, this equality still holds at termination. 
Hence $R$ is a minimal transversal. If $\mathcal F=\varnothing$, the algorithm returns the empty set. 
Since every constraint contains at most two vertices and every vertex belongs to at most two constraints, 
the running time is $O(1+|\mathcal F|)$. The array and position map for 
$\mathcal M$ allow constant-time insertion, deletion, and selection.

\begin{lemma}\label{prop:algorithm-correct}
Algorithm~\ref{alg:transversal} terminates almost surely and returns a red-stable transversal of $\mathcal D$. Each attempt succeeds with probability at least $s_0=1-17\alpha^{16}>0.9009$, where $\alpha=29/40$. The expected number of attempts is less than $1.11$.
\end{lemma}
\begin{proof}
If $\mathcal D=\varnothing$, the algorithm returns the empty set. Otherwise, \cref{lem:independent} shows that $S$ is independent in $C$ for every assignment. If $D\in\mathcal D_0$ is missed, \cref{lem:encoding} shows that $T$ avoids $Z_D$ and meets every constraint in $\mathcal F_D$. Algorithm~\ref{alg:minimal-witness} returns $R\in\mathcal R_D$ such that $E(D,R)\in\mathcal B$ is true. The main algorithm resamples exactly the variables of this event. Its inner loop is therefore a legal Moser--Tardos execution on the fixed family of distinct atomic events used in \Cref{sec:path}, without enumerating that family.

By \cref{lem:resampling-check,lem:harris-resampling}, the inner loop terminates almost surely and meets every member of $\mathcal D_0$. The proof of \cref{thm:transversal} bounds the probability of missing some member of $\mathcal D_1$ at termination by $17\alpha^{16}$. Hence the final test succeeds with probability at least $s_0$. Repeating the same event-choice rule with fresh independent samples gives independent attempts, so the expected number of attempts is at most $s_0^{-1}<1.11$. At return, the inner loop has enforced $\mathcal D_0$ and the final test has checked $\mathcal D_1$. Thus the returned set meets every member of $\mathcal D$.
\end{proof}

\begin{lemma}\label{prop:algorithm-time}
Given the cycle order and indexed vertex lists for $\mathcal D$, Algorithm~\ref{alg:transversal} uses expected time $O(n+mk)$, assuming constant expected time for an exact Bernoulli sample.
\end{lemma}
\begin{proof}
Preprocessing and initialization take $O(n+mk)$ time. 
Let $\nu$ denote the number of resamplings in one attempt. 
By the resampling-count bound and the estimate in the proof of \cref{lem:loads},
\[
\begin{aligned}
\mathbb E[\nu]
&\le \sum_{E\in\mathcal B}\mu(E)\le \sum_{D\in\mathcal D_0}
\sum_{R\in\mathcal R_D}\mu(E(D,R))\le m\alpha^k.
\end{aligned}
\]

The second inequality allows the same atomic event to be counted more than once.

Since $R$ is minimal, for each $b\in R$, there exists $F_b\in\mathcal F_D$ such that $R\cap F_b=\{b\}$. 
These sets $F_b$ are distinct, and hence $|R|\le|\mathcal F_D|$. 
Therefore,
$|U|\le |Z_D|+|R|
\le |Z_D|+|\mathcal F_D|=k$.

A resampling on $U$ changes membership in $S$ only on $U\cup N_P(U)$, 
which has size at most $3|U|$. 
Since each vertex belongs to at most $k$ members of $\mathcal D$, 
the affected counts can be updated in $O(k|U|)=O(k^2)$ time. 
Computing $R$ requires $O(k)$ time. 
Thus the expected running time of one attempt is
$O\bigl(n+mk+k^2\mathbb E[\nu]\bigr)
=O(n+mk+mk^2\alpha^k)
=O(n+mk)$,
where the last equality follows from $k\alpha^k<2/27$ by \cref{lem:loads}.

If $\mathcal D=\varnothing$, the algorithm returns immediately. Otherwise, let $T_i$ be the cost of attempt $i$, excluding preprocessing performed only once, and let $N$ be the index of the first successful attempt. Under the fixed rules and reset convention above, regard an infinite sequence of independent identically distributed attempts as defined in advance, with the algorithm stopping after attempt $N$. The event $\{N\ge i\}$ depends only on the first $i-1$ attempts and is therefore independent of $T_i$. Since the costs are nonnegative, by Tonelli\textquotesingle s theorem, we have
\begin{align*}
\mathbb E\!\left[\sum_{i=1}^{N}T_i\right]
&=\sum_{i\ge1}\mathbb E\!\left[\mathbf 1_{\{N\ge i\}}T_i\right]=\sum_{i\ge1}\Pr(N\ge i)\,\mathbb E[T_i]=\mathbb E[T_1]\,\mathbb E[N]\le s_0^{-1}\mathbb E[T_1]=O(n+mk).
\end{align*}
Here $\mathbb E[T_1]=O(n+mk)$ by the one-attempt estimate, and $\mathbb E[N]\le s_0^{-1}<1.11$ by \cref{prop:algorithm-correct}. Adding the one-time preprocessing cost proves the claimed total bound. The cost and success of a single attempt may be dependent.
\end{proof}

For an $r$-regular graph with $r\ge18$ and a specified Hamilton cycle, form the distinct closed neighborhoods in $G-E(C)$ by hashing in expected $O(nr)$ time. The reduction in \cref{lem:reduction} has $k=r-1$ and $m\le n$, so \cref{prop:algorithm-time} gives a red-stable green-dominating set in expected time $O(nr)$. 

\section*{Statement of AI}
The proofs in this paper were developed with substantial contributions from ChatGPT (5.5 and 6 Astra).
During the preparation of this paper, the authors provided it with the underlying proof framework
and key ideas, as well as relevant references, including the starting point of a refinement of Local Lemma due to Lu and Sz\'ekely \cite{LS-arxiv09} and the mixture of Thomassen's method, which results in an improvement of Haxell et al.'s result. The authors discussed with ChatGPT 5.5 to firstly improve Theorem 1.2. ChatGPT 6 suggested using the tool from \cite{HarrisSrinivasan2017} instead to further improve it to $r\geq 17$. The authors almost rewrote the paper completely and used ChatGPT (OpenAI) for language polishing, structural organization, and improving the exposition of some proofs. All mathematical arguments and proofs were independently checked and finalized by the author, who takes full responsibility for the manuscript.

\end{document}